\documentclass[a4paper]{article}
\usepackage[english]{babel}
\usepackage[utf8x]{inputenc}
\usepackage[T1]{fontenc}
\usepackage[top=2.4cm,bottom=3.2cm,left=1.8cm,right=1.8cm]{geometry}

\usepackage{graphics}
\usepackage{epsfig}
\usepackage{mathptmx}
\usepackage{times}
\usepackage{amsmath}
\usepackage{amsthm}
\usepackage{amssymb}
\usepackage{xcolor}
\usepackage{multirow}

\usepackage{cancel}
\usepackage[normalem]{ulem}

\usepackage[export]{adjustbox}
\usepackage{footnote}

\makeatletter
\let\NAT@parse\undefined
\makeatother
\usepackage{hyperref}

\def\Rset{\mathbb{R}}

\newtheorem{prop}{Proposition}

\newtheorem{assum}{Assumption}
\newtheorem{defi}{Definition}

\title{\LARGE \bf
Identifiability and observability of the SIR model with quarantine\\
}
\author{Frédéric Hamelin\thanks{F.~Hamelin is with Institut Agro, Rennes, France
		{\tt frederic.hamelin@agrocampus-ouest.fr}}, 
	Abderrahman Iggidr\thanks{A.~Iggidr is with Université de Lorraine, CNRS, Inria, IECL, F-57000 Metz, France
		{\tt Abderrahman.Iggidr@inria.fr}}, 
	Alain Rapaport\thanks{A.~Rapaport is with MISTEA, Univ.~Montpellier,
		INRAE, Institut Agro, Montpellier, France
		{\tt alain.rapaport@inrae.fr}}, 
	Gauthier Sallet\thanks{G.~Sallet is with Université de Lorraine, CNRS, IECL, F-57000 Metz, France {\tt gauthier.sallet@univ-lorraine.fr}}\,
	 and 
	Max O. Souza\thanks{M.O.~Souza is with Instituto de Matem\'atica e Estatística, Universidade Federal Fluminense, Niter\'oi - RJ, 24210-201, Brasil
	{\tt maxsouza@id.uff.br}}
}

\begin{document}

\maketitle

\begin{abstract}
We analyze the identifiability and observability of the well-known SIR epidemic model with an additional compartment Q of the sub-population of infected individuals that are placed in quarantine (SIQR model), considering that the flow of individuals placed in quarantine and the size of the quarantine population are known at any time. Then, we focus on the problem of identification of the model parameters, with the synthesis of an observer.
\end{abstract}

\section{Introduction}
Many papers in epidemiology proposing a mathematical model using dynamical systems, face the problem of parameter estimation. In general, some parameters are given, extracted from the literature, while the remaining unknown parameters are estimated by fitting the model to some observed data, usually by means of an optimization algorithm based on least squares or maximum likelihood methods.  Nevertheless, relatively few studies about the intrinsic property of a model to admit a unique set of parameters values for a given choice of measured variables. On the other hand, this a question that is well known in automatic control. Investigation of identifiability in mathematical epidemiology is relatively recent 
\cite{MR3553947,MR2726194,MR1957979,MR2817815,MR2785878,Eisenberg:2013aa,MR2142487}.
Indeed, to the best of our knowledge, the first paper considering the problem in a epidemic model studies identifiability of an intra-host model of HIV, and it has been published in 2003 in an Automatic Control journal \cite{MR1957979}. It is also surprising that the observability and identifiability of the original Kermack-Mckendrick  model has not been more studied, since this system has been widely used to model an outbreak of an infection.
The  observability and identifiability of the classical model SIR, with demography and constant population, has been first studied in 2005 \cite{MR2142487}.

The global health crisis of COVID-19 outbreak has led to a spectacular resurgence of interest in this type of models, but with specificities related to the detection and isolation of infected individuals \cite{m-y-li-SIR}. This is why we revisit here the issues of identification and observability for an extended `SIQR' model \cite{MR2003m:92074} for which such an analysis had not yet been performed (to the best of our knowledge). Once the question of identification has been settled, we also tackle the task of proposing a practical strategy for reconstructing the unique values of the parameters.

\section{The models}

Inspired by \cite{MR2003m:92074,MR2428383}, we consider the classical SIR model (see for instance \cite{MR3752193}), where $S$, $I$, $R$ denote the size of the populations of respectively susceptible, infected and recovered individuals, with an additive compartment where $Q$ denotes the size of the population of identified and isolated infectious individuals that have been removed from the infected population and placed in quarantine:
\begin{equation}
\label{model1}
\left\{\begin{array}{lll}
\dot S  & =  & -\beta S\frac{I}{N-Q} \\
\dot I & =  & \beta S\frac{I}{N-Q} - (\rho+\alpha) I\\
\dot Q  & =  & \alpha I - \rho Q\\ 
\dot R  & = & \rho I + \rho Q
\end{array}\right.\,.
\end{equation}
When the size of the total population $N$ is large and the size of the population placed in quarantine remains small compared to $N$ during the considered interval of time, one can consider a simplified model: 
\begin{equation}
\label{model2}
\left\{\begin{array}{lll}
\dot S  & =  & -\beta S\frac{I}{N} \\
\dot I & =  & \beta S\frac{I}{N} - (\rho+\alpha) I\\
\dot Q  & =  & \alpha I - \rho Q\\ 
\dot R  & = & \rho I + \rho Q
\end{array}\right.\,.
\end{equation}
Note that for both models, one has
\[
S(t)+I(t)+Q(t)+R(t)=N, \quad \forall t \geq 0\,.
\]

These models have three parameters: the infectivity parameter $\beta$, the recovery rate $\rho$, that we assume to be identical for the infected populations placed in quarantine or not, and  the rate of placement in quarantine $\alpha$. Theses parameters are unknown but we assume the following hypothesis.

\smallskip

\begin{assum} 
	\label{assum1}
	The reproduction number ${\mathcal R}_0$ verifies
	\[
	{\mathcal R}_0:=\frac{\beta}{\rho+\alpha}>1\,.
	\]
\end{assum}

\smallskip

This assumption implies that the epidemic can spread in the population i.e.~at initial time with $S(0)=N-I(0)$ close to $N$ one has $\dot I(0)>0$.


\section{The identification problem}
We assume that
\begin{itemize}
	\item the flow $\alpha I(t)$ of infected people placed in quarantine is known at any time $t\geq 0$
	\item the size $Q(t)$ of the population placed in quarantine is perfectly known at any time $t\geq 0$
	\item the size $N$ of the total population is known
	\item at initial time $0$, one has $S(0)=N-\varepsilon$, $I(0)=\varepsilon$, $Q(0)=0$, $R(0)=0$ with $\varepsilon \in (0,N)$.
\end{itemize}
We consider then the observation function
\begin{equation}
\label{observation}
y(t)=\left[\begin{array}{c} 
y_1(t)\\ 
y_2(t)
\end{array}\right]:=\left[\begin{array}{c} 
\alpha I(t)\\ 
Q(t)
\end{array}\right]
\end{equation}
and follow the usual definitions of identifiability and observability of systems \cite{MR1482525,MR1408862}.
However, note that when $Q=0$, the system is not infinitesimally identifiable: the knowledge of the outputs and all its derivative do not allow to determine formally $\rho$. At $I=0$, the system is not identifiable neither. We adopt the following definition of identifiability for these models.

\smallskip

\begin{defi}
\label{defi}
	Given $N>0$ and $\varepsilon \in (0,N)$, we shall say that system \eqref{model1} resp.~\eqref{model2} is identifiable for the observation \eqref{observation} if there exists $t>0$ such that the map
	\[
	\left[\begin{array}{c}
	\alpha\\
	\beta\\
	\rho\\
	\end{array}\right]\in \left(\Rset_{+}^\star\right)^3 \quad \longmapsto \quad y(\cdot)  \in {\mathcal C}^\infty([0,t],\Rset_+^2)
	\]
	is injective, where $(S(\cdot),I(\cdot),Q(\cdot),R(\cdot))$ is solution of the Cauchy problem for the differential system \eqref{model1} resp.~\eqref{model2} with $S(0)=N-\varepsilon$, $I(0)=\varepsilon$, $Q(0)=0$ and $R(0)=0$.
	If moreover the map
	\begin{align*}
	\left[\begin{array}{c}
	\alpha\\
	\beta\\
	\rho\\
	\varepsilon
	\end{array}\right]\in \left(\Rset_{+}^\star\right)^3\times(0,N)
	 \longmapsto \quad y(\cdot)  \in {\mathcal C}^\infty([0,t],\Rset_+^2)
	\end{align*}
	is injective, then the system \eqref{model1} resp. \eqref{model2} is identifiable and observable for the observation \eqref{observation}.

\end{defi}

\section{Analysis of the first model}

\begin{prop}
	System \eqref{model1} is identifiable and observable for the observation \eqref{observation}, in the sense of Definition \ref{defi}.
\end{prop}

\smallskip

\begin{proof}
It consists in showing that parameters and unmeasured variables $S$ and $I$ can be expressed as functions of the successive derivatives of the output vector $y$. As the variable $I$ cannot reach $0$ in finite time, we shall assume $I\neq 0$ in the following.

\smallskip

Note first that with $Q(0)=0$ one has $\dot Q(0)>0$ and then $y_2(t)=Q(t)>0$ for any $t>0$. The dynamics of $Q$ gives directly the expression of the parameter $\rho$ as:
\begin{equation}
\label{rho}
\rho= \, \dfrac{ y_1(t)-\dot{y}_2(t)}{y_2(t)}\,, \quad t>0 .
\end{equation}
Posit $h_1:=  \dfrac{\dot y_1 }{y_1}$. One has from the dynamics of $I$
\begin{equation}\label{eq1}
h_1 = \dfrac{\beta \, S}{N-Q}-\alpha- \rho  .
\end{equation}
and then
\begin{equation}\label{eq2}
(N-Q) \, \dot{h}_1=-\dfrac{\beta^2\, S}{N-Q} \, I + \dfrac{\beta \, S}{N-Q}\, \dot Q.
\end{equation}
Using the equality $\dfrac{\beta \, S}{N-Q}=\alpha + h_1+\rho$  from \eqref{eq1}, one obtains
\begin{equation}\label{eq3}
h_2:= (N-y_2) \, \dot{h}_1 = (h_1+\alpha +\rho)\, ( -\beta \,I+  \dot Q).
\end{equation}
Let us write the derivative of $h_2$:
\begin{align*}
& \dot h_2 = \dot h_1\, (-\beta \, I + \dot Q)\\ 
& \quad + (h_1+\alpha+\rho) \, \left [   -\beta \,  \dfrac{\beta \, S}{N-Q}\, I + \beta \, (\alpha+\rho) \, I  + \ddot Q\right ]
\end{align*}
which can be also expressed as
\begin{align*}
& \dot h_2 = \dot h_1\, (-\beta \, I + \dot Q)\\ 
& \qquad + (h_1+\alpha+\rho) \, \left [   -\beta \, I  \, (h_1+\alpha+\rho)+ \beta \, (\alpha+\rho) \, I  + \ddot Q\right ] \\
& \quad = \dot h_1\, (-\beta \, I + \dot Q) + (h_1+\alpha+\rho)  \, \left [   -h_1\, \beta \, I   + \ddot Q\right ] .
\end{align*}
Then, using relation \eqref{eq3}, one obtains the expression
\[
\dot h_2= \dot h_1 \,  (-\beta \, I + \dot Q)+ \dfrac{h_2}{(-\beta\, I +\dot Q)}\, \left [   h_1\, (-\beta \, I  + \dot Q)  - h_1\, \dot Q +\ddot Q\right ]
\]
which implies
\begin{align*}
& (-\beta \, I  + \dot Q) \,  \dot h_2 = \dot h_1 \,  (-\beta \, I + \dot Q)^2\\
& \qquad + h_2\, h_1 \,  (-\beta \, I  + \dot Q)+ h_2\, ( - h_1\, \dot Q +\ddot Q)
\end{align*}
or equivalently the equation
\begin{equation*}
  \dot h_1 \,  (-\beta \, I + \dot Q)^2+ (h_2\, h_1 -\dot h_2)  (-\beta \, I  + \dot Q)+ h_2\, ( - h_1\, \dot Q +\ddot Q) =0
\end{equation*}
to be fulfilled.

\medskip

Observe that this last equation is a second order polynomial in the variable  $X=-\beta \, I + \dot Q$.
From \eqref{eq2} and $\mathcal R_0>1$ one has
\begin{equation}
\label{ineq1}
\left\{\dot h_1\right\}_{t=0} =\dfrac{\beta \, \varepsilon }{N }\, (-\beta +\alpha)\left(1-\frac{\varepsilon}{N}\right) <0
\end{equation}
and this allows us to show that one also has
\begin{equation}
\label{ineq2}
\left\{h_2\, ( - h_1\, \dot Q +\ddot Q)\right\}_{t=0} =\alpha \, \beta \, \rho \, \varepsilon^2 (\beta -\alpha)\left(1-\frac{\varepsilon}{N}\right) >0
\end{equation}
Indeed, one has 
\begin{align*}
&\ddot Q = \alpha \beta \dfrac{S\, I}{N-Q}-\alpha \, (\rho+\alpha)\, I -\rho \, \dot Q\\
& \qquad \Rightarrow \left\{\ddot Q \right\}_{t=0}= \alpha\, \varepsilon \, \left(\beta\left(1-\frac{\varepsilon}{N}\right)-\alpha-2 \rho\right)
\end{align*}
With $h_1(0) =\beta\left(1-\frac{\varepsilon}{N}\right)-\rho-\alpha$, one obtains
\[
\begin{array}{lll}
\left\{( -h_1\, \dot Q + \ddot Q)\right\}_{t=0} & = & \alpha\, \varepsilon \, (\beta\left(1-\frac{\varepsilon}{N}\right)-\alpha-2 \rho)\\
& & \qquad - (\beta\left(1-\frac{\varepsilon}{N}\right) -\rho-\alpha) \, \alpha \, \varepsilon\\
& = & -\alpha  \rho \,\varepsilon <0 
\end{array}
\]
and with $h_2(0)= \beta \,\varepsilon\,(-\beta + \alpha)\left(1-\frac{\varepsilon}{N}\right) $, one gets
\[ 
\left\{\left(h_2\, ( -h_1\, \dot Q + \ddot Q)\,\right)\right\}_{t=0} =  \alpha \, \beta \, \rho\, \varepsilon^2 (\beta -\alpha)\left(1-\frac{\varepsilon}{N}\right) >0 .
\]
Observe also that one has $X(0)=-\beta \, I(0) + \dot Q(t)=\varepsilon\, (-\beta + \alpha) <0$. Therefore, by continuity w.r.t.~$t$, we obtain that for $t>0$ small enough, $X$ is the unique negative solution of 
\[
\dot h_1 \,  X^2+ (h_2\, h_1 -\dot h_2)\,  X + h_2\, ( - h_1\, \dot Q +\ddot Q) =0 .
\] 
that is
\[
X=\frac{-2(h_2\, h_1 -\dot h_2)-\sqrt{(h_2\, h_1 -\dot h_2)^2+4\dot h_1 ( h_1\, \dot y_2 -\ddot y_2)}}{2\dot h_1}\,.
\]
The parameter $\alpha$ can be then obtained from equation \eqref{eq2}
\[
\alpha= \frac{(N-y_2)\dot h_1}{X}-h_1-\rho
\]
where $\rho$ is given by \eqref{rho}. The initial condition $\varepsilon$ is simply reconstructed by $\varepsilon=y_1(0)/\alpha$ and finally one obtains the parameter $\beta=\alpha-X(0)/\varepsilon$.
\end{proof}

\section{Analysis of the simplified model}

\begin{prop}
	System \eqref{model2} is identifiable and observable for the observation \eqref{observation}, in the sense of Definition \ref{defi}.
\end{prop}

\smallskip

\begin{proof}
	As for model \eqref{model1}, one can determine the parameter $\rho$ from any positive time as
	\[
	\rho = \frac{y_1(t)-\dot y_2(t)}{y_2(t)} , \quad t > 0 .
	\]
	Then from the dynamics of $I$ one can write
	\begin{equation}
	\label{h1}
	\frac{\beta S(t)}{N}-\alpha = h_1(t):=\frac{\dot y_1(t)}{y_1(t)} + \rho , \quad t >  0
	\end{equation}
	where $h_1$ is a known function. Differentiating $h_1$ with respect to the time gives
	\begin{equation}
	\label{h1dot}
	\dot h_1(t)=-\beta^2S(t)\frac{I(t)}{N^2}=-\frac{\beta I(t)}{N}(h_1(t)+\alpha)
	\end{equation}
	and differentiating twice
	\[
	\ddot h_1(t)= -\frac{\beta}{N}\left(\frac{\beta S(t)}{N}-\rho-\alpha\right)I(t)(h_1(t)+\alpha)-\frac{\beta I(t)}{N}\dot h_1(t) .
	\]
	With the expression \eqref{h1}, we rewrite this last equation as follows
	\[
	\ddot h_1(t)= -\frac{\beta I(t)}{N}(h_1(t)-\rho)(h_1(t)+\alpha)-\frac{\beta I(t)}{N}\dot h_1(t)
	\]
	and with expression \eqref{h1dot}
	\[
	\ddot h_1(t)= -\dot h_1(t)(h_1(t)-\rho)-\frac{\beta I(t)}{N}\dot h_1(t) .
	\]
	One obtains then the following expression for $\beta I(t)$
	\[
	\beta I(t)=-N\left(\frac{\ddot h_1(t)}{\dot h_1(t)} +h_1(t)-\rho\right)
	\]
	(note that this expression is well defined because $h_1(t)+\alpha  > 0$ for any $t$ and thus $\dot h_1(t)<0$).
	
	\medskip
	
	Finally, from \eqref{h1dot} one reconstructs the parameter
	\[
	\alpha=-\frac{N\dot h_1(t)}{\beta I(t)}-h_1(t)
	\]
	and then the parameter
	\[
	\beta = \alpha\frac{\beta I(t)}{y_1(t)}
	\]
	At last, the initial condition is recovered as $\varepsilon=y_1(0)/\alpha$.
\end{proof}

\section{Parameter estimation}

The former analyses have shown that models \eqref{model1} and \eqref{model2} are not infinitesimally identifiable at time $0$ when the initial state is $(N-\epsilon,\epsilon,0,0)$. One has to wait a short time $t>0$ to have $Q(t)>0$ and formally identify parameters. Thus, we expect a weak accuracy of the parameters estimation at the very beginning, that should improve with time while the state get away from this initial state and new measurements come. This is why we have opted for a dynamical estimation with the help of observers. The use of observers, although usually dedicated to state estimation (rather to parameters estimation) possesses the advantage to tune the speed of error decay. Moreover, the choice of a speed-accuracy compromise can be balanced thru simulations with synthetic data corrupted by noise, before facing real data. 

Note also that for large times, the solutions of \eqref{model1} and \eqref{model2} converge asymptotically to non-observable non-identifiable states when $I$ and $Q$ are both null. Consequently, we do not look precisely for results about asymptotic convergence of the error (as this is usually done in observers theory), but rather for an exponential decay of the error estimation during initial transients.

In this section, we shall consider the additional hypothesis

\smallskip

\begin{assum} 
	\label{assum2}
	One has
	\[
	\alpha \leq \rho
	\]
\end{assum}
that means that the rate of placement in quarantine is not larger than the recovery rate, which is often the case in epidemic regimes.

\medskip

We shall denote the elementary symmetric polynomials, where $X$ is a vector in $\Rset^n$, as
\[
\sigma_k^n(X):=\sum_{1\leq i_1< \cdots < i_k\leq n}\left( \prod_{j=1}^k X_j\right), \quad i=1\cdots n
\]
\begin{prop} Let $\lambda$ and $\mu$ be two positive vectors in $\Rset^4$ and $\Rset^3$ respectively. For $t>0$, consider the dynamical system
\begin{equation}
\label{full_observer}
\left\{\begin{array}{l}
    \dot {\hat {z}}_1 = \hat \delta -\hat\rho
    -K_1(\hat{z_{1}} -\log(y_1(t)),\\[2mm]
   \dot {\hat {z}}_2  =  \dfrac{y_1(t)}{y_2(t)} -\hat\rho -K_2(\hat{z_{1}} -\log(y_1(t)),\\[2mm]
    \dot {\hat \delta}= -K_3(\hat{z_{1}} -\log(y_1(t))), 
    \\[2mm]
    \dot {\hat \rho}= -K_4(\hat{z_{1}} -\log(y_1(t)))-(\hat{z_{2}} -\log(y_2(t))) ,  \\[2mm]
    \dot {\hat{y}}_1 =  \hat{v} y_1(t) -K_5 y_1(t)(\hat{y_1}-y_1(t)), \\[2mm]
 \dot {\hat{v}} = -\hat{k} \dfrac{y_1(t)}{N}-K_6 y_1(t)(\hat{y_1}-y_1(t)),\\[2mm]
  \dot {\hat{k}} = - K_7Ny_1(t)(\hat{y_1}-y_1(t)) 
\end{array}\right.
\end{equation}
with the gains vector
\[
K=\left[\begin{array}{c}
\sigma_1^4(\lambda)\\
\sigma_1^4(\lambda)+\sigma_3^4(\lambda)\\
-\sigma_4^4(\lambda)\\ 
-\sigma_2^4(\lambda)-\sigma_4^4(\lambda)-1\\
\sigma_1^3(\mu)\\
\sigma_2^3(\mu)\\
-\sigma_3^3(\mu)
\end{array}\right]
\]
Then, the output vector
\begin{equation}
    \label{estimator}
\left[\begin{array}{c}
\hat\rho(t)\\[2mm]
\hat\beta(t) = \dfrac{1}{2}\left(\hat k(t)-\sqrt{\max(\hat k(t)^2-4\hat\delta(t)\hat k(t),0)}\right)\\[2mm]
\hat \alpha(t)=\hat\beta(t)-\hat\delta(t)
\end{array}\right]
\end{equation}
is an estimator of $[\rho,\beta,\alpha]^\top$, whose exponential decay of the error can be made as fast as desired keeping the number
\[
l=\min_{i,j} \min(\lambda_i,\mu_j) .
\]
large, as long as $S(t)/N$ remains close to $1$. Moreover, the state $I$ is estimated with
\begin{equation}
    \label{Iestim}
\hat I(t) = \frac{y_1(t)}{\hat\alpha(t)}
\end{equation}
\end{prop}

\medskip

\begin{proof}
Posit
\begin{equation}
    \label{defdelta}
\delta:=\beta-\alpha .
\end{equation}
As far as $S/N$ remains close to $1$, the size of the population $I$ is small compared to $S$, and the dynamics of the outputs $y_1=\alpha I$ and $y_2=Q$ can be approximated by the linear dynamics
\begin{equation}
\label{ApproxLin1}
\left\{\begin{array}{l}
    \dot y_1 =  \delta\,y_1- \rho y_1 \\[2mm]
    \dot y_2=   y_1-\rho \, y_2\\
\end{array}\right.
\end{equation}
For $t>0$ and $i=1,2$, consider the new outputs $z_i(t)=\log(y_i(t))$, whose dynamics is given by the system
\begin{equation}
\label{ApproxLin2}
\left\{\begin{array}{l}
    \dot z_1 =\delta- \rho  \\[2mm]
    \dot z_2=  \exp{(z_1-z_2)}-\rho \\
\end{array}\right.
\end{equation}
For this sub-system with unknown parameters $\delta$ and $\rho$, we consider the following candidate observer in $\Rset^4$:
\begin{equation}
\label{obs1}
\left\{\begin{array}{l}
    \dot {\hat {z}}_1 = \hat \delta -\hat\rho
    -K_1(\hat{z_{1}} -z_{1}),\\[2mm]
   \dot {\hat {z}}_2  =  \exp{(z_1-z_2)} -\hat\rho -K_2(\hat{z_{1}} -z_{1}),\\[2mm]
    \dot {\hat \delta}= -K_3(\hat{z_{1}} -z_{1})\\[2mm]
    \dot {\hat \rho}= -K_4(\hat{z_{1}}-z_1) -(\hat{z_{2}} -z_{2}) ,   
\end{array}\right.
\end{equation}
The dynamics of the error $e_1=(\hat{z_{1}},\hat{z_{2}},\hat\delta,\hat\rho)^\top-(z_1,z_2,\delta,\rho)^\top$ 
is given by the linear time invariant system $\dot e_1 =M_1\,e_1$ with 
\[M_1= 
\left(\begin {array}{cccc} 
-K_1& 0 & 1 & -1
\\[3mm] 
-K_2& 0 &0 & -1\\[3mm] 
-K_3 & 0 & 0 & 0\\[3mm] 
-K_4& -1 &  0 & 0
\end {array} 
\right) 
\]
A calculation shows that with the choice
\[\left\{ 
\begin{array}{l}
K_{1}=\lambda_1+\lambda_2+\lambda_3+\lambda_4, \\[2mm]
K_{2}=\sum\lambda_i
+\sum\lambda_i\lambda_j\lambda_k 
,\\[2mm]
K_{3}=-\lambda_1\,\lambda_2\,\lambda_3\,\lambda_4
{K_{42}}\\[2mm]
K_{4}=-\left(\sum\lambda_i\lambda_j
+ \lambda_1\lambda_2\lambda_3\lambda_4+1
\right),
\end{array}\right.
\]
the spectrum of $M_1$ is $\{-\lambda_i, \; i =1\cdots 4\}$. This shows that the first four equations of system \eqref{full_observer} gives the reconstruction of the parameters $\delta$ and $\rho$ with an exponential decay of the error larger than $\min_i \lambda_i$.

\medskip

Posit now
\begin{equation}
 \label{defk}
 k:=\frac{\beta^2}{\alpha} .
\end{equation}
and consider the variable
\begin{equation}
 \label{defv}
 v(t):=\beta\frac{S(t)}{N}-\rho-\alpha
 \end{equation}
As far as $S(T)/N$ remains close to $1$, one can write the approximation
 \[
 \dot v = -\frac{\beta^2}{\alpha}\frac{S(t)}{N^2}y_1(t) \simeq -\frac{\beta^2}{\alpha N}y_1(t)
 \]
Then, this amounts to approximate the dynamics of \eqref{model1} or \eqref{model2} in the $(y_1,v,k)$ coordinates by the following dynamical system
\[
\left\{\begin{array}{l}
\dot y_1 = v y_1\\[2mm]
\dot v = -\dfrac{k}{N}y_1
\end{array}\right.
\]
where $k$ is an unknown parameter, for which we consider the following candidate observer in $\Rset^3$
\begin{equation}
\label{obs2}
\left\{
\begin{array}{l}
 \dot {\hat{y}}_1  =  \hat{v} y_1 -K_5 y_1 (\hat{y_1}-y_1) \\[3mm]
 \dot {\hat{v}}  = -\hat{k} \dfrac{y_1}{N}-K_6 y_1 (\hat{y_1}-y_1)\\[3mm]
  \dot {\hat{k}}  = -K_7N y_1(\hat{y_1}-y_1) 
\end{array}
\right.
\end{equation}
whose dynamics of the error $e_2=(\hat{y_1},\hat{v},\hat{k})^\top-(y_1,v,k)^\top$ is given by 
the non-autonomous linear system
\begin{equation}
    \label{dynerror2}
    \dot e_2 = y_1(t) M_2\, e_2
\end{equation}
with
\[
M_2=\left(\begin{array}{ccc}
-K_5 & 1 & 0 \\[3mm]
-K_6 & 0 & -\dfrac{1}{N} \\[3mm]
-K_7N & 0 & 0
\end{array}\right)\]
One can easily check that for the choice
\[\left\{
\begin{array}{l}
 K_5 =  \mu_1 + \mu_2 + \mu_3, \\[3mm]
 K_6 = \mu_1\mu_2 + \mu_1\mu_3 + \mu_2\mu_3,\\[3mm]
 K_7 = -\mu_1\mu_2\mu_3 ,
\end{array}
\right.
\]
the spectrum of $M_2$ is $\{-\mu_i , \; i =1\cdots 3\}$. Then, from \eqref{dynerror2}, we obtain the upper bound on the error decrease
\[
|\hat k(t)-k| \leq \exp\left(-(\min_j \mu_j) \int_0^T y_1(\tau)d\tau\right)||e_2(0)|| 
\]
whose exponential decay can be made as large as desired with large $l$.

\medskip

Finally, from the reconstruction of parameters $\delta$, $k$ by observers \eqref{obs1}, \eqref{obs2} and expressions \eqref{defdelta} and \eqref{defk}, the original parameters $\alpha$, $\beta$ are recovered as roots of 
\begin{equation}
    \label{roots}
\beta^2-k\beta+k\delta=0 \Rightarrow \beta = \frac{k\pm\sqrt{k^2-4k\delta}}{2}
\end{equation}
Note first that Assumption 1 implies $\beta>\alpha$ and thus $k^2-4k\delta>0$. Moreover, one has $k>\beta(1+\frac{\rho}{\alpha})$, and by 
Assumption 2 one has $k>2\beta$, which implies that only the smaller root of \eqref{roots} is valid, leading to the expression \eqref{estimator} of the estimator. 
Note that this expression preserves the exponential decay of the error obtained for $\delta$, $k$ and $\rho$.
\end{proof}

\medskip

Let us make some comments about this observer. It consists in reconstructing functions of the parameters $\delta$ and $k$ and not directly the parameters $\alpha$, $\beta$. There is an apparent redundancy of variables $\hat z_1$ and $\hat y_1$ in dynamics \eqref{full_observer}, which reconstruct $\log y_1$ and $y_1$. Indeed, this allows to decouple the observer into two sub-systems of dimensions $4$ and $3$, which avoids the use of two large correction gains compared to a full order observer. Finally, outputs of these two sub-systems are coupled in expression \eqref{estimator} to reconstruct the original parameters.

\section{Numerical illustrations}
\begin{figure}[h!]
	\begin{center}
		\includegraphics[width=7cm]{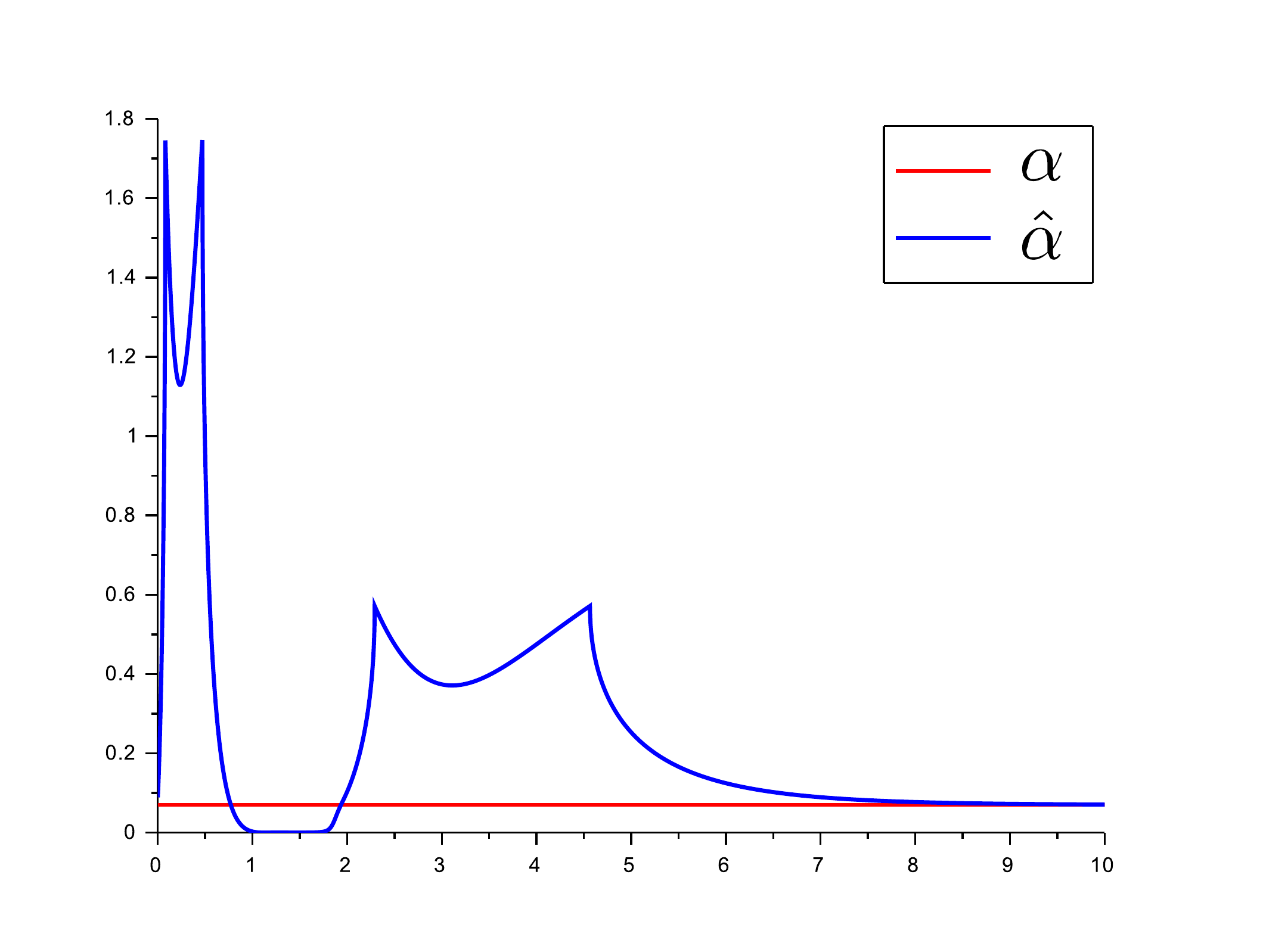}
		\includegraphics[width=7cm]{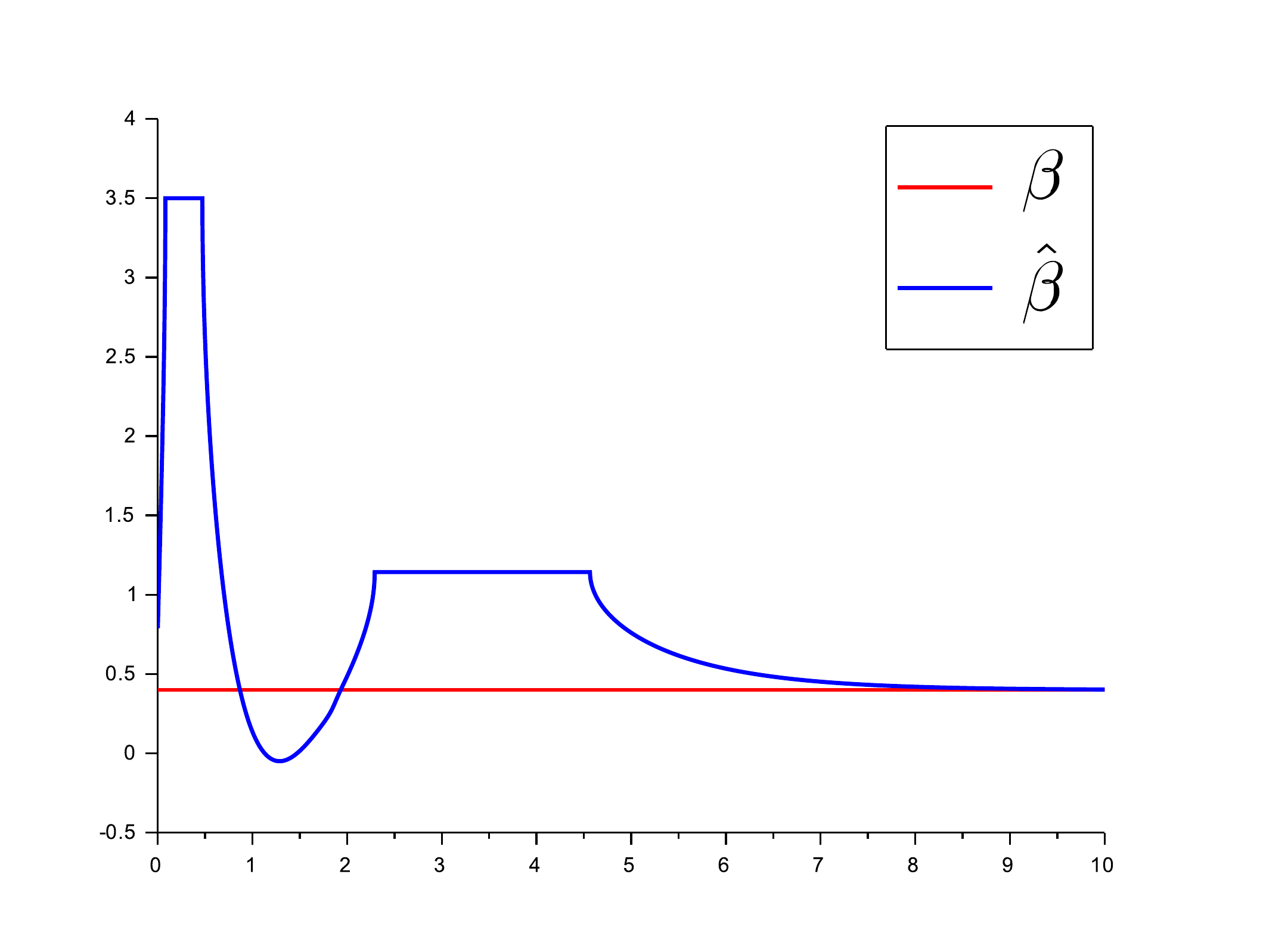}
		\includegraphics[width=7cm]{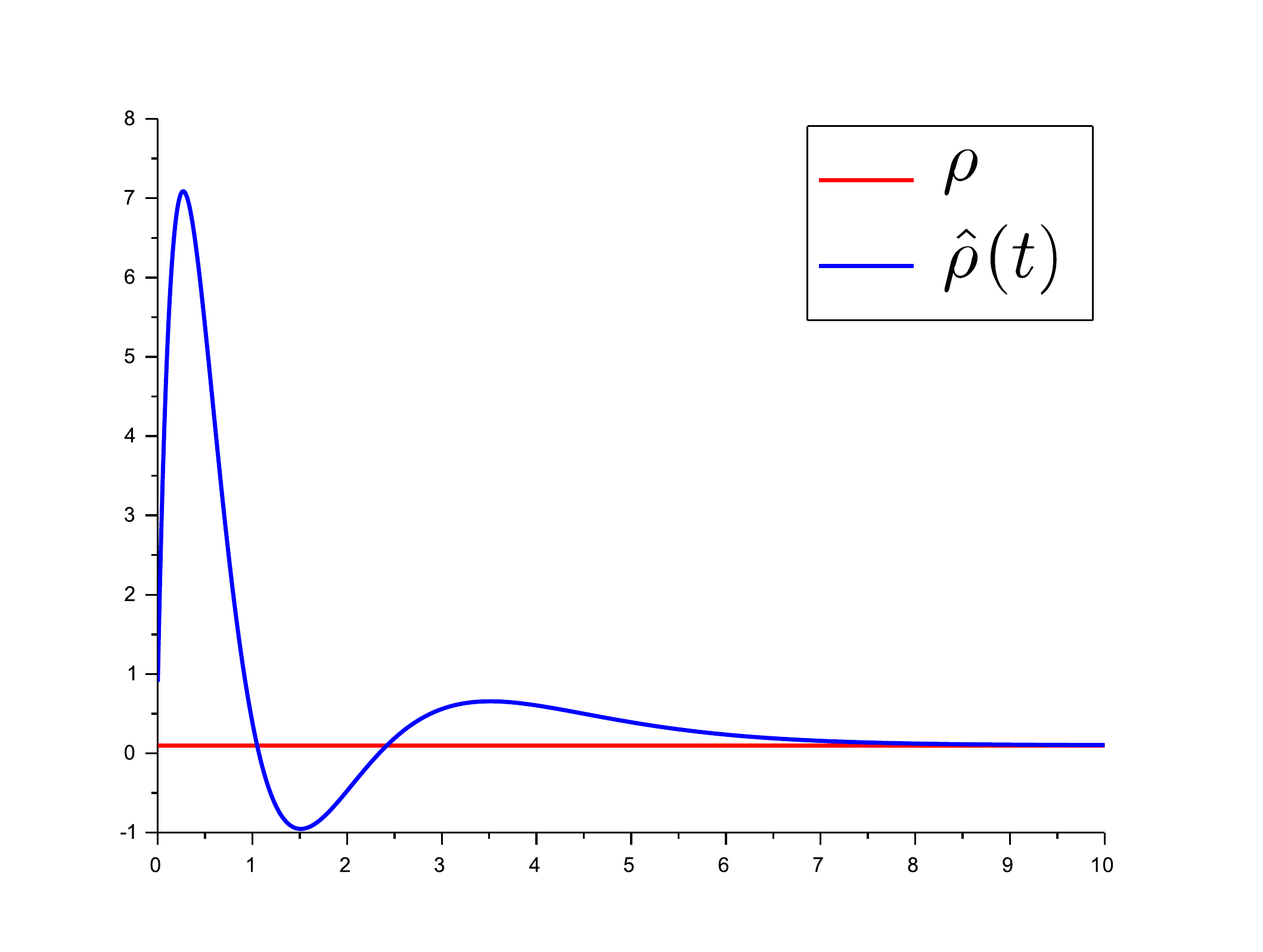}
		\includegraphics[width=7cm]{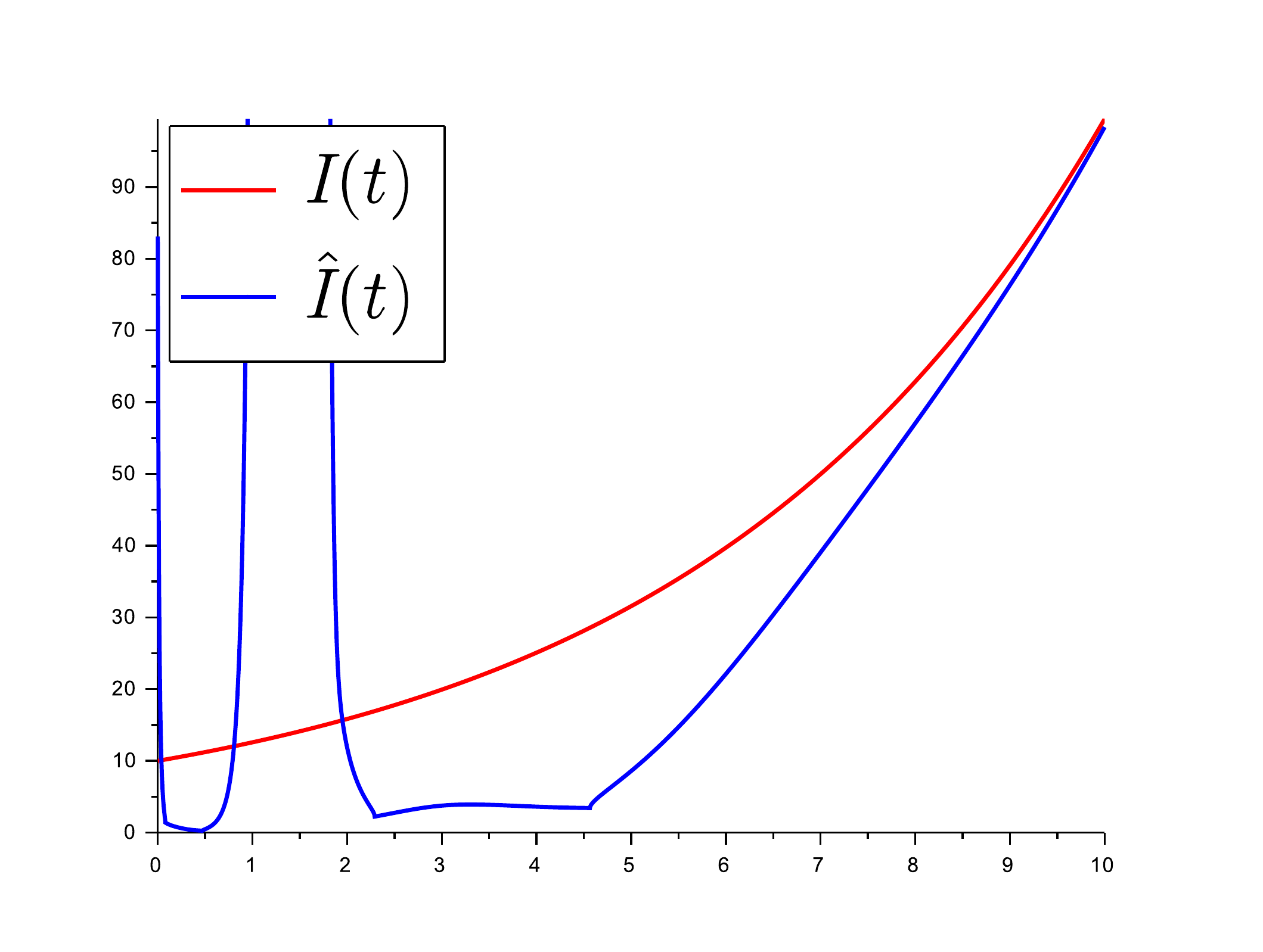}
		\caption{Observer simulation without measurement noise}
		\label{fig-sansbruit}
	\end{center}
\end{figure}

The proposed observer has been tested with synthetic data
for a population size $N=10^5$ with parameter values $\alpha=0.07$, $\beta=0.4$, $\rho=0.1$ and initial condition
$I(0)=10$, $Q(0)=5$, $R(0)=0$ over a time horizon of $10$ days (see Figure \ref{fig-sansbruit}). The gains have been computed for
the choice of vectors $\lambda=[1;1.5;2;2.5]$ and $\mu=[1/(13.10^3);1/(15.10^3);1/(19.10^3)]$. Note that vector $\mu$ has been chosen quite small to avoid too large gains when multiplied by $N$ in the observer equations.

\begin{figure}[h!]
	\begin{center}
		\includegraphics[width=7cm]{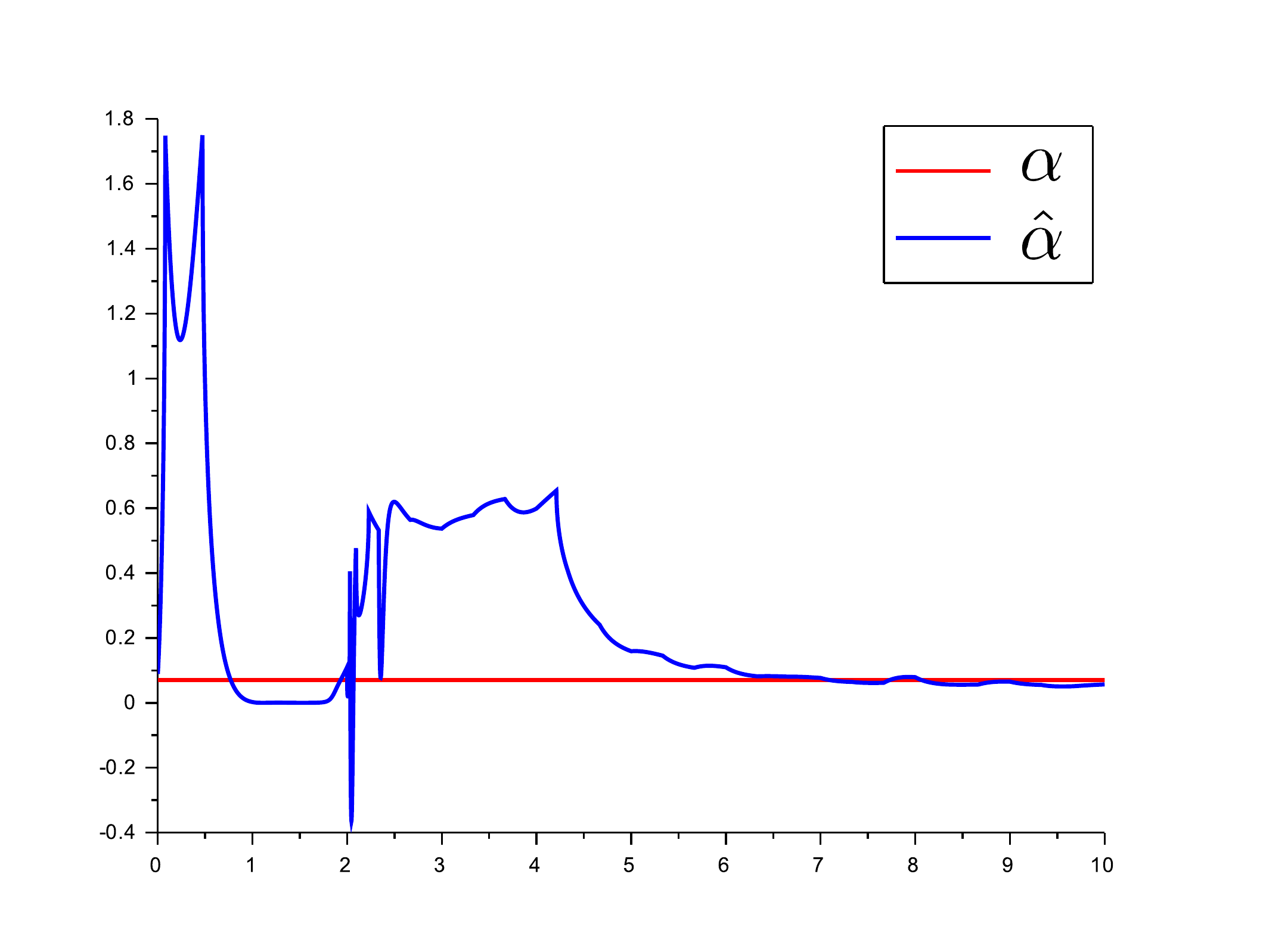}
		\includegraphics[width=7cm]{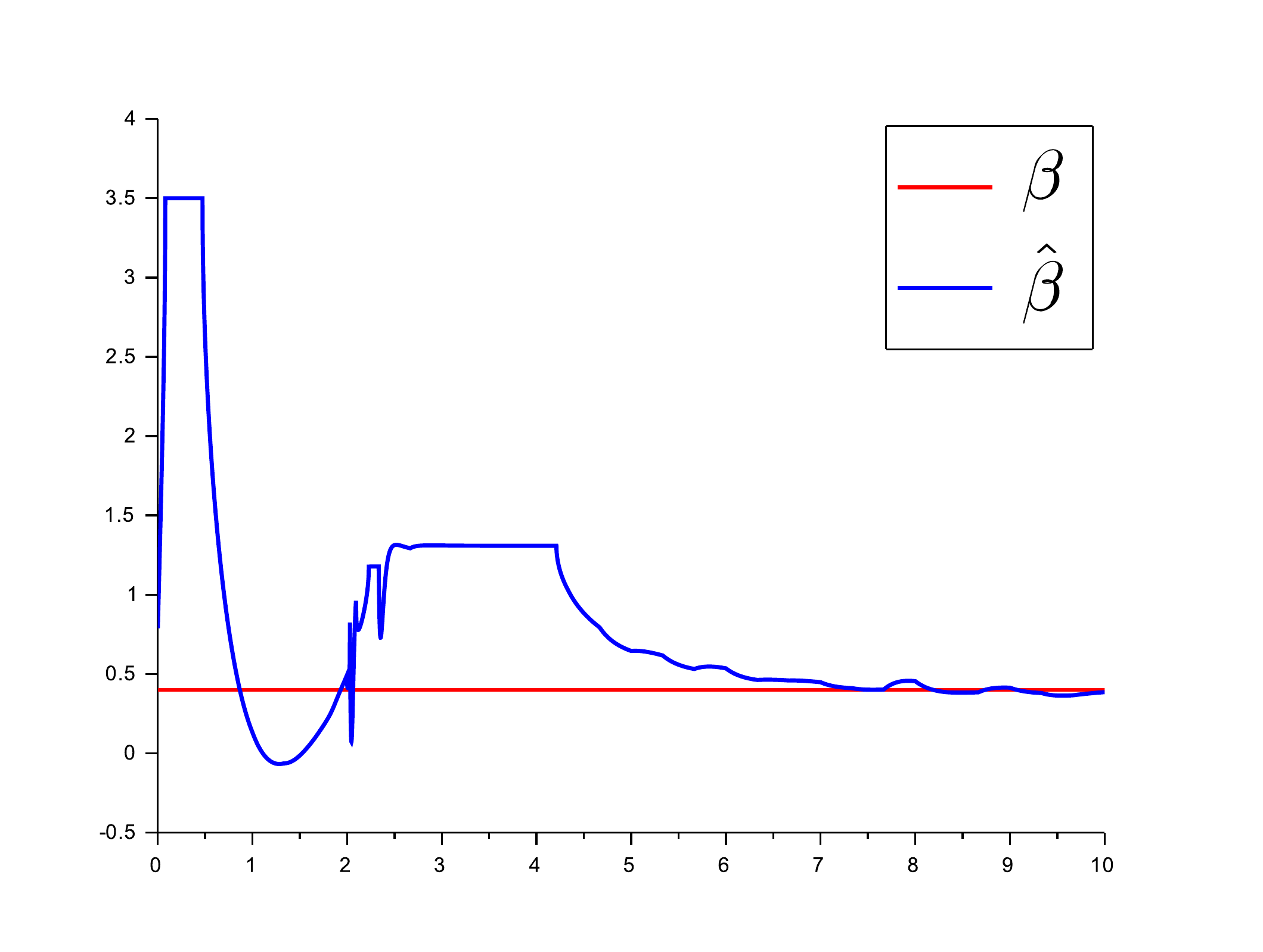}
		\includegraphics[width=7cm]{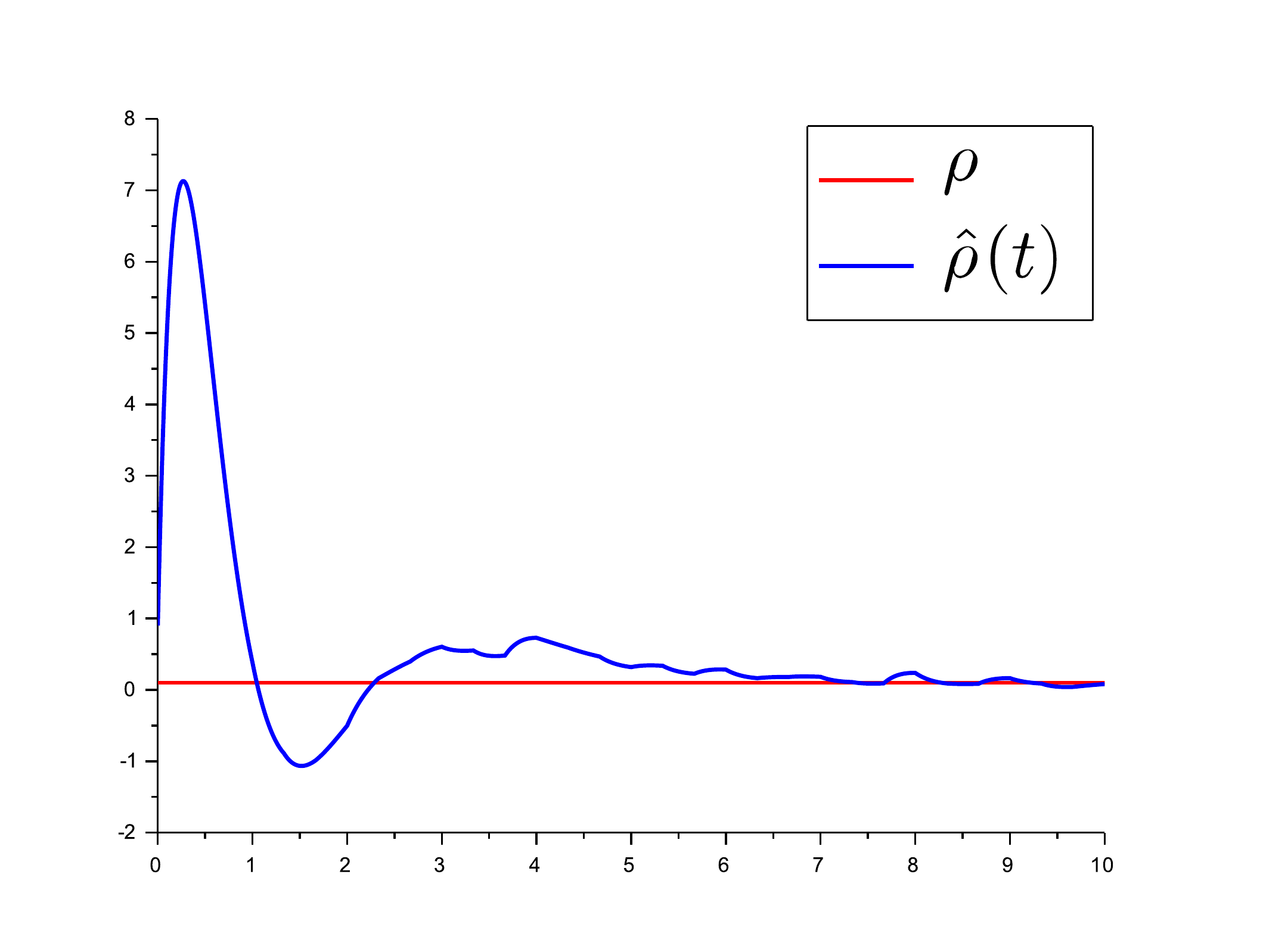}
		\includegraphics[width=7cm]{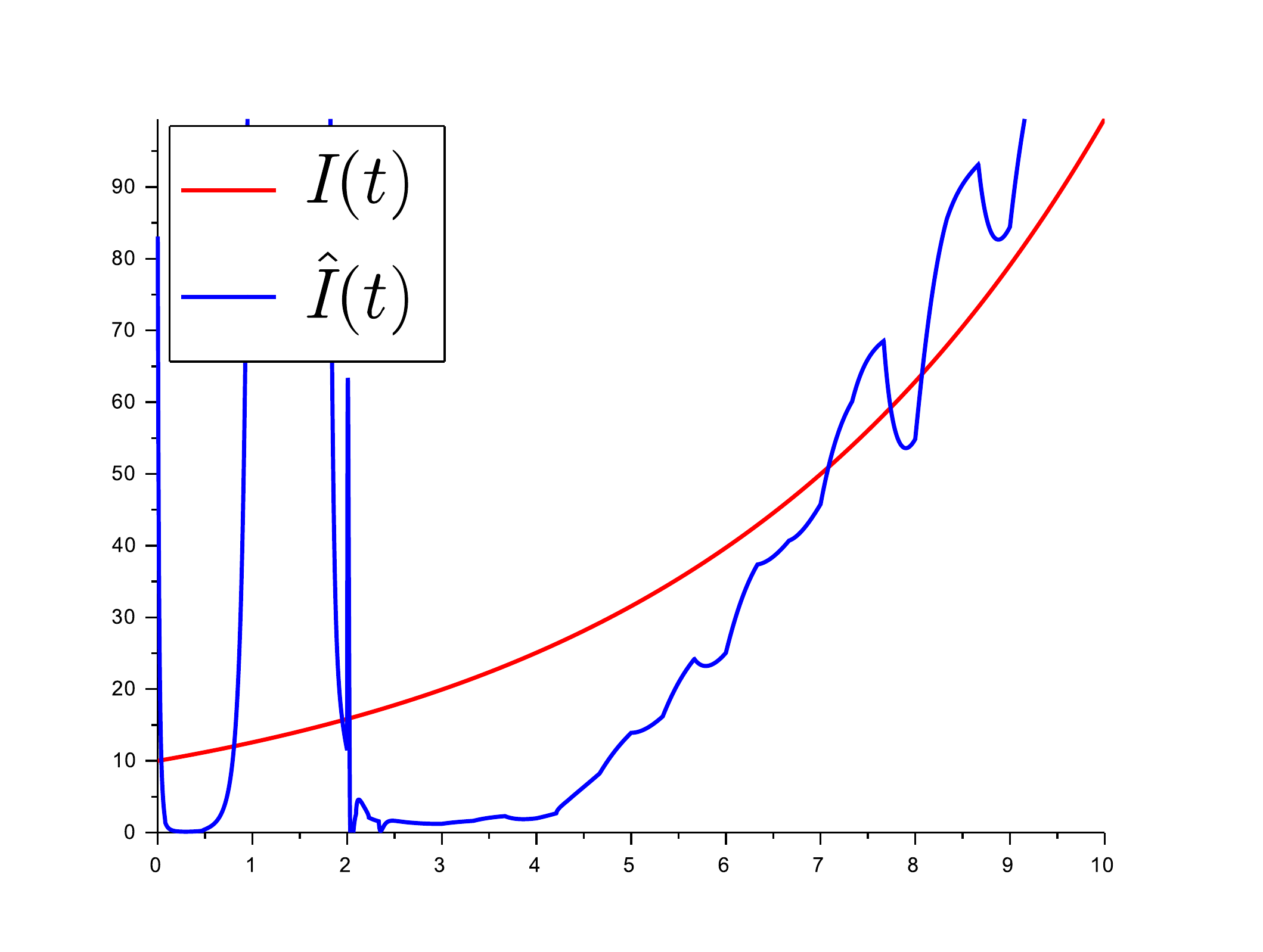}
		\caption{Observer simulation with measurement noise}
		\label{fig-avecbruit}
	\end{center}
\end{figure}

\begin{figure}[h!]
	\begin{center}
		\includegraphics[width=6cm]{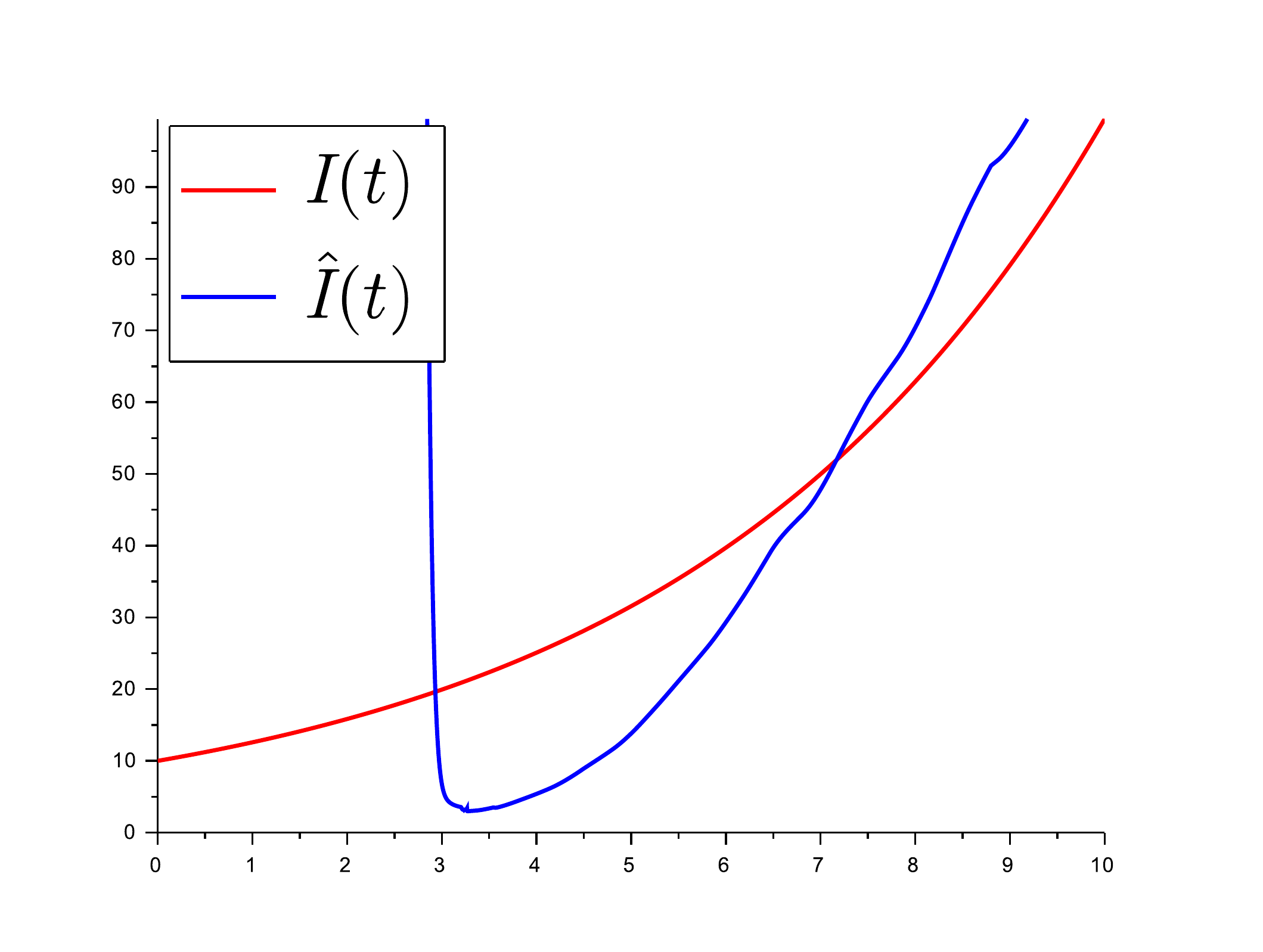}
		\caption{Smoothed estimation of the infected population in presence of measurement noise}
\label{figIlisse}
	\end{center}
\end{figure}

Then, we have simulated a measurement noise with a centered Gaussian law of variance equal to $5\%$ of the signal (see Figure \ref{fig-avecbruit}).
Because the expression \eqref{Iestim} of the estimation of the state $I$ is not filtered, we have applied a moving average smoothing to the estimation $\hat I$ (see Figure \ref{figIlisse}).

Finally, these simulations show that the method allows to reconstruct the parameter values in few days in a quite accurate manner. The estimation of the size of the infected population $I$ allows then the use of the model for predictions of the epidemics.

\section{Conclusion}
This work shows that although the identifiability of the SIR-Q models has singularity points where measured variables are null, it is possible to design an observer with exponential decay of the estimation error during the first stage of the epidemics, and recover parameters in few days. Further investigations will concern real data of COVID epidemics provided by various territories.

\clearpage

\bibliographystyle{abbrv}
\bibliography{observer-min.bib}

\end{document}